\PassOptionsToPackage{hyphens}{url}	
\documentclass[]{siamonline171218}
\usepackage{amsfonts}
\usepackage{graphicx}
\usepackage{epstopdf}
\ifpdf%
  \DeclareGraphicsExtensions{.eps,.pdf,.png,.jpg}
\else
  \DeclareGraphicsExtensions{.eps}
\fi
\usepackage{amsopn}

\newtheorem{example}[theorem]{Example}

\newcommand{\TheTitle}{%
  A low-rank matrix equation method for solving PDE-constrained optimization problems\thanks{This version dated May 11, 2020.}
}

\newcommand{\ShortTitle}{%
  Low-rank matrix equation method for PDE-constrained optimization
}

\newcommand{\TheAuthors}{Alexandra B\"unger, Valeria Simoncini, and Martin Stoll}

\headers{\ShortTitle}{\TheAuthors}
\title{{\TheTitle}}

\author{Alexandra B\"unger\thanks{Technische Universit\"at Chemnitz, \email{alexandra.buenger@mathematik.tu-chemnitz.de}}
\and 
Valeria Simoncini\thanks{Alma Mater Studiorum - Universit\`a di Bologna, Italy, and
IMATI-CNR, Pavia. \email{valeria.simoncini@unibo.it}}
\and
Martin Stoll\thanks{Technische Universit\"at Chemnitz, \email{martin.stoll@mathematik.tu-chemnitz.de}}
}
\ifpdf
\hypersetup{
	pdftitle={\TheTitle},
	pdfauthor={\TheAuthors}
}
\fi

\usepackage{amsopn}
\usepackage{booktabs}
\usepackage{algorithmicx}
\usepackage{amsmath}
\usepackage{subcaption}

\usepackage{booktabs}

\newcommand{\R}{\mathbb{R}}
\newcommand{\eps}{\varepsilon}
\newcommand{\D}{\,\mathrm{d}}

\def\minres{{{\sc Minres }}}
\def\newmethod{{{\sc sys2mateq }}}
\def\lrminres{{{\sc lrminres }}}
\def\matlab{{{\sc matlab }}}
\def\intel{{{\sc Intel }}}

\begin{document}

\maketitle


\begin{abstract}
PDE-constrained optimization problems arise in a broad number of applications such as hyperthermia cancer treatment or blood flow simulation. Discretization of the optimization problem and using a Lagrangian approach result in a large-scale saddle-point system, which is challenging to solve, and acquiring a full space-time solution is often infeasible. We present a new framework to efficiently compute a low-rank approximation to the solution by reformulating the KKT system into a Sylvester-like matrix equation. This matrix equation is subsequently projected onto a small subspace via an iterative rational Krylov method and we obtain a reduced problem by imposing a Galerkin condition on its residual. In our work we discuss implementation details and dependence on the various problem parameters. Numerical experiments illustrate the performance of the new strategy also when compared to other low-rank approaches.
\end{abstract}

\begin{keywords}
  PDE-constrained Optimization, Matrix Equation, Rational Krlyov Subspace
\end{keywords}

\begin{AMS}
	65F10, 65F50, 15A69, 93C20 
\end{AMS}

\section{Introduction}\label{sec:intro}
A vast number of todays technological advancements is due to the ever advancing improvement in modeling and solving of partial differential equation (PDE) problems. One class of especially challenging problems is PDE-constrained optimization which has a broad number of applications, reaching from engineering design and control problems to medical applications like cancer treatment \cite{schenk2009,biegler2003large,fisher2009data}.

We consider a typical PDE-constrained optimization problem on a space-time cylinder $[0,T] \times \Omega$ given by
\begin{align} 
\min_{y,u} J(y,u)  \label{eq:objective0} \\
\mbox{subject to } \dot{y} = \mathcal{L}(y) + u, \label{eq:constraint0}
\end{align}
where we are interested in the minimization of the functional $J(y,u)$ constrained by a differential operator $\mathcal{L}(y)$ with space and time dependent state $y$ and control $u$. An extensive analysis of this type of problems can be found e.g. in \cite{Troeltzsch2010} or \cite{Hinze2009}. In this work, we will follow the popular approach of first discretizing the problem and then formulating the discrete first-order optimality conditions as described in \cite{book::IK08}. Finding an efficient numerical solution to the large system of equations resulting from these conditions is of high interest and many approaches exist to solve the resulting saddle point system such as using a block preconditioner to subsequently solve the preconditioned system with iterative solvers like \minres (cf. \cite{stoll2010,pearson2012regularization,rees2010optimal,schoberl2007symmetric}).

However, todays PDE models often result in millions of variables and solving optimization problems governed by such huge models is still a challenging task due to their size and complexity. They often even prove to be impossible to solve on a normal computer due to the memory required to store the usually large and dense solution computed by standard algorithms. Therefore, memory efficient solvers, which compute reduced approximations to the solution are a desirable approach to 
solve these challenging problems. We here propose a widely applicable low-rank solver, 
which computes a projection of the solution onto a small subspace. 

The method relies on the reformulation of the problem into a matrix equation, which will be illustrated in section~\ref{sec:pb}. 
The numerical solution of this PDE-constrained optimization problem poses a significant challenge regarding the complexity of both storage and computational resources. We show that it is possible to separate spatial and temporal data via a novel low-rank matrix equation solver in section~\ref{sec:lowrank}. We introduce a method to reduce the system of equations into a 
single low-rank matrix equation. 
{In section~\ref{section:comp_specifications1} we discuss the choice of the
approximation space, while in section~\ref{section:comp_specifications2} we analyze} stopping 
criteria and a fast update scheme of the residual computation. And finally in section~\ref{sec:num}, we examine the performance and applicability of our method on various numerical examples. We show the method's robustness with respect to the various problem parameters as well as assess the performance of our method compared with a previously proposed low-rank approach on a distributed heat equation model as well as a boundary control problem. Further, we show the applicability of our method to more challenging problems such as a partial state observation and a non-symmetric PDE constraint.

\section{Problem formulation} \label{sec:pb}

We consider a PDE-constrained optimal control problem on a space-time cylinder with time interval $[0,T]$ and a spatial domain $\Omega$ with a subdomain $\Omega_{1} \subseteq \Omega$ as in Equation \eqref{eq:objective0}-\eqref{eq:constraint0}. The functional we want to minimize reads,
	\begin{equation}
	J(y,u) =  \frac{1}{2} \int_{0}^T \int_{\Omega_1} (y - \hat{y} )^2 \D x \D t + \frac{\beta}{2} \int_{0}^T \int_{\Omega} u^2 \D x \D t, \label{eq:objective}
	\end{equation}
	where $y$ is the state, $\hat{y}$ the given desired state, and $u$ the control, which is regularized by the control cost parameter $\beta$. For the PDE subject to which we want to minimize the functional $J(y,u)$ let us exemplarily consider the heat equation with $\mathcal{L}(y) = \Delta y$ and Dirichlet boundary condition,
	\begin{align}
	\dot{y} -\Delta y &= u \mbox{ in } \Omega, \label{eq:PDE}\\
	y &= 0 \mbox{ on } \partial \Omega. 
	\end{align}
	There are two distinctive options to proceed with this problem. Either we formulate the optimality conditions and discretize them subsequently or we first discretize the problem and then formulate the discrete first-order optimality conditions, known as the Karush-Kuhn-Tucker (KKT) conditions \cite{book::IK08}. In this work we follow the second approach.
	
	We discretize space and time in an all-at-once approach. The spatial discretization is done with finite elements and to discretize in time, we split the time interval into $n_T$ intervals of length $\tau = \frac{T}{n_T}$. Using a rectangle rule, the discretization of \eqref{eq:objective} becomes
	\begin{equation}
	 \sum_{t=1}^{n_T}  \frac{\tau}{2} (y_{t} - \hat{y}_{t})^T M_1 (y_{t} - \hat{y}_t) + \frac{\tau \beta}{2} u_t^T M u_t,
	\end{equation}
	where $y_t, \hat{y}_t$, and $u_t$ are spatial discretizations of $y,\hat{y}$, and $u$ of size $n$ for each time step $t = 1,\hdots,n_T$. Using an implicit Euler-scheme the discrete formulation of the PDE \eqref{eq:PDE} reads
	\begin{align}
	 \frac{M (y_{t} - y_{t-1})}{\tau} + K y_t &= M u_t, \mbox{ for } t = 1,\hdots, n_T,\\
	 y_0 &= 0.
	\end{align}
    Here, $M,N,M_1 \in \R^{n\times n}$ are the mass matrices of the spatial discretization and $K\in \R^{n \times n}$ denotes the so-called stiffness matrix resulting from the discretization of the differential operator $\mathcal{L}(y)$. The boundary constraints are incorporated into the stiffness matrix.  On further details about the discretization of PDE operators we refer the reader to \cite{elman2005}.
    We collect the discretizations of the variables in matrices $Y = [y_1,\hdots,y_{n_T}] \in \R^{n \times n_T}$ and  denote their vectorizations by $\underline{Y}= \mbox{vec}(Y)$, respectively for $\hat{y}$ and $u$. With this 
we write the optimization problem in compact form as
	\begin{align}
	\min_{Y,U} \frac{\tau}{2} (\underline{Y}-\underline{\hat{Y}})^T \mathcal{M}_1 (\underline{Y}-\underline{\hat{Y}}) &+ \frac{\tau \beta}{2} \underline{U}^T \mathcal{M} \underline{U}, \label{eq:oc_discrete1} \\
	\mbox{s. t. } \mathcal{K}\underline{Y} - \tau \mathcal{N} \underline{U} &= 0, \label{eq:oc_discrete2}
	\end{align}
	with the matrices
	\begin{align}
	 \mathcal{M} = \begin{bmatrix} M & & & \\ & M & & \\ & & \ddots & \\ & & & M \end{bmatrix},\qquad \mathcal{K} = \begin{bmatrix} M + \tau K & & & \\ -M & \ddots & & \\ & \ddots & \ddots \\ & & -M & M + \tau K \end{bmatrix},
	\end{align}
	where mass matrices $M\in \R^{n\times n}$ and stiffness matrix $K \in \R^{n\times n}$ repeat $n_T$ times each. The matrices $\mathcal{N}$ and $\mathcal{M}_1$ are block diagonal matrices like $\mathcal{M}$ but with different mass matrices $N$ and $M_1$ respectively on the diagonal.
	 
	 To solve the discrete problem \eqref{eq:oc_discrete1} - \eqref{eq:oc_discrete2}, we have to solve the system of equations resulting from the first order optimality conditions \cite{benzi_2005}, They state that an optimal solution must be a saddle point of the Lagrangian of the problem,
	 \begin{equation}
	  \nabla \mathcal{L}(\underline{Y}^*, \underline{U}^*, \underline{\Lambda}^*) = 0.
	 \end{equation}
    The Lagrangian of this problem reads
	\begin{equation}
	\mathcal{L}(\underline{Y},\underline{U},\underline{\Lambda}) = \frac{\tau}{2}(\underline{Y}-\underline{\hat{Y}})^T \mathcal{M}_1 (\underline{Y}-\underline{\hat{Y}}) + \frac{\tau \beta}{2} \underline{U}^T \mathcal{M} \underline{U} + \underline{\Lambda}^T (\mathcal{K} \underline{Y} - \tau \mathcal{N} \underline{U}).
	\end{equation}
	Thus, the optimal solution solves the following set of linear equations,
	\begin{align} 
	0 = \nabla_Y \mathcal{L}(\underline{Y},\underline{U},\underline{\Lambda)} &= \tau \mathcal{M}_1 (\underline{Y} - \underline{\hat{Y}}) + \mathcal{K}^T \underline{\Lambda}, \label{eq:KKT_conditions1} \\
	0 = \nabla_U \mathcal{L}(\underline{Y},\underline{U},\underline{\Lambda)} &= \tau \beta \mathcal{M} \underline{U} - \tau \mathcal{N}^T \underline{\Lambda}, \label{eq:KKT_conditions2} \\
	0 = \nabla_\Lambda \mathcal{L}(\underline{Y},\underline{U},\underline{\Lambda}) &= \mathcal{K} \underline{Y} - \tau \mathcal{N} \underline{U}. \label{eq:KKT_conditions3}
	\end{align}
	Memory-wise, an effective approach to solve this large-scale problem is to use a low-rank approach, which finds a cheap representation of the solution matrix
    within a low-rank subspace range($V$) and a reduced solution $Z$, whose number of columns is related to the number of employed time steps,
    \begin{equation}
     Y \approx V Z ,
    \end{equation}
 and similarly for the other matrix variables $U$ and $\Lambda$.
In this work we show that we can compute such a low-rank solution more efficiently than recent approaches by rearranging the large system of equations into a generalized Sylvester matrix equation of the form
	 \begin{equation} \label{eq:matrix_equation}
    A_1 X + X C + A_2 X I_0 - A_3 X D - F_1F_2^T = 0.
    \end{equation}
    The resulting matrix equation can be efficiently solved by using a tailored low-rank Krylov-subspace method. This new approach greatly reduces storage and time requirements while being robust to parameter changes.
This matrix equation oriented methodology allows one to clearly identify space and time dimensions, and aims at reducing the problem dimension in the space variables; the time variable is handled using an all-at-once procedure. Similar strategies have been used, for instance, in \cite{Breiten.Simoncini.Stoll.16,Palitta.19}.
    
    One assumption we make 
is that we have a low-rank approximation 
or representation of the desired state as
    \begin{equation}
     \hat{Y} \approx Y_1 Y_2^T,
    \end{equation}
    with $Y_1 \in \R^{n \times r}$, $Y_2 \in \R^{n_T \times r}$ and $r < n_T$. First, we introduce
 the auxiliary matrices
    \begin{align}
	 \mathcal{I} = \begin{bmatrix} 1 & & \\ & \ddots & \\ & & 1 \end{bmatrix} \mbox{ and } C = \begin{bmatrix} 1 & & & \\ -1 & 1 & & \\ & \ddots & \ddots & \\ & & -1 & 1 \end{bmatrix},
	\end{align}
    to rewrite the system matrices as Kronecker products $\mathcal{M} = \mathcal{I} \otimes M$, $\mathcal{M}_1 = \mathcal{I} \otimes M_1$, $\mathcal{N} = \mathcal{I} \otimes N$ and $\mathcal{K} = \mathcal{I} \otimes \tau K + \mathcal{C} \otimes M$. With this our system of KKT conditions \eqref{eq:KKT_conditions1} - \eqref{eq:KKT_conditions3} becomes
    \begin{align}
    \tau (\mathcal{I} \otimes M_1) \underline{Y} + (\mathcal{I} \otimes \tau K^T + \mathcal{C}^T \otimes M^T) \underline{\Lambda} - \tau (\mathcal{I} \otimes M_1) \underline{\hat{Y}} &=0, \label{eq:KKT_kronecker1}\\
    \tau \beta (\mathcal{I} \otimes M) \underline{U} - \tau (\mathcal{I} \otimes N^T) \underline{\Lambda} &= 0, \label{eq:KKT_kronecker2}\\
    (\mathcal{I} \otimes \tau K + \mathcal{C} \otimes M) \underline{Y} - \tau (\mathcal{I} \otimes N) \underline{U} & = 0. \label{eq:KKT_kronecker3}
    \end{align}
	We exploit the relation
	\begin{equation} \label{eq:kronecker_relation}
     (W^T \otimes V) \mbox{vec}(Y) = \mbox{vec}(VYW)
    \end{equation}
    to rewrite the equations \eqref{eq:KKT_kronecker1} - \eqref{eq:KKT_kronecker3} as
    \begin{align}
     \tau M_1 Y + \tau K^T \Lambda + M^T\Lambda C - \tau M_1 \hat{Y} &= 0, \label{eq:system1} \\
     \tau \beta M U - \tau N^T \Lambda &= 0, \label{eq:betau} \\
     \tau K Y + M Y C^T - \tau N U &= 0. \label{eq:system3}
    \end{align}
    The mass matrix $M$ arising from a standard Galerkin method will always be symmetric and can be considered lumped, i.e., diagonal. Therefore, we can eliminate Equation \eqref{eq:betau} by setting  
 $U = \frac{1}{\beta} M^{-1}N^T\Lambda$ in the remaining two equations,
    \begin{align}
     \tau M_1 Y + \tau K^T \Lambda + M^T \Lambda C - \tau M_1 \hat{Y} &= 0, \label{eq:orig1}\\
     \tau K Y + M Y C^T - \frac{\tau}{\beta} N M^{-1} N^T \Lambda &= 0,  \label{eq:orig2}
    \end{align}
    which are equivalent to
    \begin{align}
     M^{-1} M_1 Y + M^{-1}K^T \Lambda +  \Lambda \tilde{C} - M^{-1}M_1 \hat{Y} & = 0, \\
     M^{-1} K Y + Y \tilde{C}^T - \frac{1}{\beta} M^{-1}N M^{-1}N^T \Lambda &= 0,
    \end{align}
    where $\tilde{C} = \frac{1}{\tau}C$. For now, let us assume that $K = K^T$ but our approach can be easily generalized for non-symmetric $K$ as we will demonstrate later on. Now this representation corresponds to
    \begin{equation}
    \begin{split}
     M^{-1}K \begin{bmatrix} Y & \Lambda\end{bmatrix} + \begin{bmatrix} Y & \Lambda\end{bmatrix} \underbrace{\begin{bmatrix} \tilde{C}^T & 0 \\ 0 & \tilde{C} \end{bmatrix}}_{C_1} + M^{-1}M_1 \begin{bmatrix} Y & \Lambda\end{bmatrix} \underbrace{\begin{bmatrix} 0 & \mathcal{I} \\ 0 & 0 \end{bmatrix}}_{I_0} &+ \\
     M^{-1}N M^{-1}N^T \begin{bmatrix} Y & \Lambda\end{bmatrix} \underbrace{\begin{bmatrix} 0 & 0 \\ \frac{-\mathcal{I}}{\beta} & 0 \end{bmatrix}}_{D} - \underbrace{\begin{bmatrix} 0 & M^{-1}M_1 \hat{Y} \end{bmatrix}}_{F_1 F_2^T} &= 0.
    \end{split}
    \end{equation}
    We denote $X = \begin{bmatrix} Y & \Lambda\end{bmatrix}$, $A_1 = M^{-1}K$, $A_2 = M^{-1}M_1$, $A_3 = M^{-1}N M^{-1}N^T$, $F_1 = M^{-1}M_1 Y_1$ and $F_2 = [0_{n_T \times r}, Y_2]$, where $\hat Y = Y_1 Y_2$ with $Y_1, Y_2$ of low column and row rank, respectively. Then we get the desired format from Equation \eqref{eq:matrix_equation} as
    \begin{align}
     A_1 X + X C_1 + A_2 X I_0 + A_3XD - F_1 F_2^T = 0, \label{eq:matrix_equation_final}
    \end{align}
    where the left-hand coefficient matrices have size $n\times n$ while the right-hand ones have size $2n_T\times 2n_T$, so that $X \in {\mathbb R}^{n\times 2n_T}$.

    \section{Low rank solution}\label{sec:lowrank}
The generalized Sylvester equation in (\ref{eq:matrix_equation_final}) replaces the
large Kronecker product based system of equations \eqref{eq:KKT_kronecker1} - \eqref{eq:KKT_kronecker3}.
Since the solution matrix $X \in \R^{n\times 2n_T}$ will be dense and potentially very large, to exploit the new setting it is desirable to find an appropriate approximation space and a low-rank reduced matrix approximation $Z\in \R^{p \times 2n_T}$ such that
    \begin{equation} \label{eq:lowRank_X}
     X \approx V_p Z.
    \end{equation}
where the orthonormal columns of $V_p \in \R^{n \times p}$ generate the approximation space.
    With this setting we can construct a reduced version of the 
matrix equation \eqref{eq:matrix_equation_final}. Let us assume we compute an approximation as in \eqref{eq:lowRank_X}.
Then the residual matrix associated with (\ref{eq:matrix_equation_final})
reads
    \begin{align} \label{eq:residual}
    R = A_1 V_p Z + V_p Z C_1 + A_2 V_p Z I_0 + A_3 V_p Z D - F_1 F_2^T.
    \end{align}
We impose the Galerkin orthogonality of the residual matrix with respect to the approximation space, which in the matrix inner product is equivalent to writing $V_p^T R = 0$, so that our equation becomes
    \begin{equation}
    V_p^T A_1 V_p Z + V_p^T V_p Z E + V_p^T A_2 V_p Z I_0 - V_p^T F_1 F_2^T = 0.
    \end{equation}
 Let us denote the reduced $p\times p$ coefficient matrices as
$A_{1,r} := V_p^T A_1 V_p$, $A_{2,r} := V_p^T A_2 V_p$, $A_{3,r} := V_p^T A_3 V_p$ and set $F_{1,r} = V_p^T F_1 \in \R^{p \times 1}$. The resulting reduced equation
    \begin{equation} \label{matrix_equation_reduced}
    A_{1,r} Z + \mathcal{I}_p Z C_1 + A_{2,r} Z I_0 + A_{3,r} Z D - F_{1,r} F_2^T = 0
    \end{equation}
    has the same structure as the original matrix equation \eqref{eq:matrix_equation_final} but its size is reduced to $p \times 2n_T$. By exploiting once again the relation in Equation \eqref{eq:kronecker_relation}, we get the small linear system of equations
    \begin{equation} \label{eq:lineq}
     \big ( (\mathcal{I}_{2n_T} \otimes A_{1,r}) + (C_1^T \otimes \mathcal{I}_p)  + (I_0^T \otimes A_{2,r}) + (D^T \otimes A_{3,r}) \big ) \underline{Z} = F_{1,r}F_2^T,
    \end{equation}
    with $\underline{Z} = \mbox{vec}(Z)$.
    For a small subspace size $p \ll n$ this system of equations is significantly easier to solve and we can either use a direct or an iterative method to do so.
 If the obtained approximate solution $V_pZ$ is not sufficiently good, then the space can be expanded and a new approximation constructed, giving rise to an iterative method. {The use of low-rank methods 
within optimization of large-scale systems has been successfully documented in several articles,
and we refer the reader to \cite{stoll2014,dolgov2016fast,cichocki2017tensor,benner2016block} for recent accounts.}

    \section{Subspace computation} \label{section:comp_specifications1}
    To construct the projection \eqref{eq:lowRank_X} we need an iterative subspace method, which constructs a relevant subspace for our problem. For this we make use of rational Krylov subspaces
    \begin{equation} \label{eq:Krylov}
    V_p(A,v,s) = \mbox{span} \big \{ v, (A+s_1I)^{-1}v, \hdots, \prod_{j=1}^{p-1}(A+s_j I)^{-1}v \big \}
    \end{equation}
    with shifts $s_j$ as described in \cite{Simoncini2016}.  This approach has proven to be very well suited to solve similarly structured problems as the shifts allow for efficient updates of the subspace using relevant information on the matrices' eigenvalues.
    As an initial vector we take the right-hand side, $v_0 = F_1$, and to 
construct the subspace \eqref{eq:Krylov} we employ a tailored strategy to adapt
to the different settings. More precisely:

i) {\it Case $M_1=M$, $N$ square, full rank.}  This corresponds to a setup where desired state and control are both distributed equally on the whole domain $\Omega$. We use the matrix $A = A_1$. We observe that in this case $A_2=I$ and $A_3$ is a diagonal nonsingular matrix.

ii) {\it Case $M_1\ne M$, $N$ square, full rank.} This corresponds to a setup where, e.g., the resulting state 
is only observed on a partial domain. We construct a mixed subspace where we add the following
two new vectors in step $k$,
    \begin{equation} \label{eq:subspace}
    \{ (A_1 + s^{(1)}_k I)^{-1}v_k, \,\, (A_2 + s^{(2)}_k I)^{-1}v_k \},
    \end{equation}
so that the space dimension grows by at most two per iteration, instead of one.

iii) {\it Case $M_1= M$, $N$ tall.} In this case, we can only control a partial domain, e.g., 
the boundary. Here, $A_2=I$ while $A_3=M^{-1} N M_b^{-1} N^T$ is not invertible.
We thus define $A_3(\alpha) = A_3 + \alpha A_1$, this selection is justified below.
In this case we use the following
two new vectors in step $k$,
    \begin{equation} \label{eq:subspace3}
     \{  (A_1 + s^{(1)}_k A_3(\alpha))^{-1}v_k, \,\, ( A_3(\alpha) + s^{(2)}_k A_1)^{-1}v_k \},
    \end{equation}
and once again the space dimension grows by at most two per iteration, instead of one.

iv) {\it Case $M=M_1$, $N$ square, full rank, $K \neq K^T$.} For a number of PDE constraints, like convection-diffusion problems, the stiffness matrix is non-symmetric. In this case, we have to slightly modify the first component $A_1X$ 
 in \eqref{eq:matrix_equation_final} to $A_1 X \begin{bmatrix} \mathcal{I} & 0 \\ 0 & 0 \end{bmatrix} + M^{-1}K^TX \begin{bmatrix} 0 & 0 \\ 0 & \mathcal{I} \end{bmatrix}$. Here, again we construct a mixed subspace where we add two new vectors in step $k$, 
\begin{equation}
 \{ (A_1 + s_k^{(1)}I)^{-1}v_k, \,\, (M^{-1}K^T + s_k^{(2)}I)^{-1}v_k \}.
\end{equation}

A strategy similar to (ii)-(iv) was successfully adopted in \cite{Powell2017} for a multiterm linear matrix equation in 
a different context. 
The effectiveness of the mixed approaches in (ii)-(iii) relies on the fact that the generated space be rich in eigencomponents of both matrices. In (ii) where $A_2$ is diagonal and singular with a bounded and tight nonzero spectral interval,
a good spectral homogeneity in the space construction can be reached by shifting the matrix $A_2$ by some
$\alpha$ so that the spectra of $A_1$  and $A_2 + \alpha_2 I_n$ are contained by the same interval.
Hence, we introduce a shift parameter $\alpha_2$  in (\ref{eq:matrix_equation_final})
to get
    \begin{equation} \label{eq:shift}
    A_1 X + X (C_1 - \alpha_2 I_0) + (A_2 + \alpha_2 I_n) X I_0 - A_3 X D - F_1 F_2^T = 0.
    \end{equation}
With these premises, a good choice for $\alpha_2$ is a rough approximation to
the largest eigenvalue of $A_1$. Due to the good properties of the transformed
spectrum, the shifts were computed by only using the projection of the matrix $A_1$ onto the generated space; more details on the shift selection can be found in \cite{simoncini2011}. In the cases i) and iv), $A_2$ is not singular but still the spectra of $A_1$ and $A_2$ differ greatly. Applying the same strategy to shift $A_2$ as in \ref{eq:shift} also proved to be beneficial in these cases.

In (iii), the structure and singularity of $A_3$ was significantly more challenging, because $\|A_3\|$ inversely depends on $\beta$, hence the spectrum of a shifted version of $A_3$ may not be enclosed in that of $A_1$.
We propose to consider the equivalent problem 
   \begin{equation}
    A_1 X(I-\alpha_3 D) + X C_1 + A_2 X I_0 - (A_3+\alpha_3 A_1) X D - F_1 F_2^T = 0.
    \end{equation}
where $\alpha_3 = \frac 1 {\sqrt{\beta}}\|A_3\|_F/\|A_1\|_F$. With this formulation, we found particularly successful the use of the projection of the pair $(A_1, A_3+\alpha_3 A_1)$ onto the current approximation
space to determine the next shifts during the iteration.
    
{\it Further subspace reduction.}    
    By computing $X \approx V_p Z$ we want to reduce storage requirements and computation time. 
This goal is only achieved if the approximation space dimension remains considerably smaller than $n_T$. 
By enriching the subspace as outlined above, however,
the space dimension may in principle grow up to $n$ in case the solution is not sufficiently accurate.
To keep the subspace as small as possible when it is not optimal we include a truncation scheme which bounds the subspace dimension to less than $n_T$. 
Thus, in the worst case the solution will have maximum rank. 
 To reduce the dimension we compute a singular value decomposition 
of $Z = U \Sigma V^T$, with $\Sigma={\rm diag}(\sigma_1, \ldots, \sigma_p)$
 where  $\sigma_1 \geq \hdots, \geq \sigma_p$ are the singular values of $Z$. 
In case some of the singular values are
too small, say $\sigma_i < \eps$ for some $i$, it is clear that 
some of the subspace information is not needed for computing a good approximation to the solution. 
If this occurs we truncate the subspace with $\tilde{p} = i-1$ by setting
    \begin{equation}
    V_{\tilde{p}} = V_p U_{\tilde{p}},
    \end{equation}
    with $U_{\tilde{p}}$ denoting the first $\tilde{p}$ columns of $U$.
    Note that the next subspace additions \eqref{eq:subspace} are still conducted 
using $v_k$ from the latest generated vectors in the original subspace $V_p$, and orthogonalized with respect to the reduced subspace.
We refer the reader to the  discussion of Table~\ref{table:memory} for an illustration
of the achievable memory savings by adopting this strategy.

    \section{Residuals and stopping criteria} \label{section:comp_specifications2}
To determine a computable stopping criterion we monitor the residuals of the two equations \eqref{eq:orig1} and \eqref{eq:orig2},
   \begin{align}
   R_1 &= \tau M_1 Y + \tau K^T \Lambda + M^T \Lambda C - \tau M_1 \hat{Y}, \label{eq:res1} \\
   R_2 &= \tau K Y + M Y C^T - \frac{\tau}{\beta} N M^{-1}N^T \Lambda, \label{eq:res2}
   \end{align}
   which are closely related to the original system. 
   
   Computing the residuals poses a bottleneck in this scheme as straightforward computation is time consuming and would require forming the full solution $[Y \Lambda] = V_p Z$. To avoid this and substantially speed up the computation time, we rely on a low-rank representation of the residual and calculate its norm making use of an updated QR method. The matrix residuals can be rewritten as
   \begin{align}
   R_1 &= \tau M_1 V_p Z_Y + \tau K^T V_p Z_{\Lambda} + M V_p Z_{\Lambda} C - \tau M_1 \hat{Y} \\
   & = \begin{bmatrix} \tau M_1 V_p & \tau K^T V_p & M V_p & -\tau M_1 Y_1 \end{bmatrix} \cdot \begin{bmatrix} Z_Y & Z_{\Lambda} & Z_{\Lambda} C & Y_2 \end{bmatrix}, \\
   R_2 &= \tau K V_p Z_Y + M V_p Z_Y C^T - \frac{\tau}{\beta} N M^{-1} N^T V_p Z_{\Lambda}  \\
   & \begin{bmatrix} \tau K V_p & M V_p & - \frac{\tau}{\beta} N M^{-1} N^T V_p \end{bmatrix} \cdot \begin{bmatrix} Z_Y & Z_Y C^T & - Z_{\Lambda} \end{bmatrix}^T,
   \end{align}
   where $Z_Y$ denotes the $Y$-component of $Z$ and $Z_{\Lambda}$ denotes the component associated with $\Lambda$ respectively.
   Let any of the two quantities be denoted by 
   \begin{equation}
   R = R_L R_R^T.
   \end{equation}
   Here the matrices have dimensions $R_L \in \R^{n \times (3p+1)}$ and $R_R \in \R^{n_T\times (3p+1)}$. We consider a reduced QR decomposition of the tall and skinny matrix $R_L$,
   \begin{equation}
   R_L  = Q_1 R_1, \quad {\rm with} \quad
   Q_1 \in \R^{n \times (4p+1)},\quad R_1 \in \R^{(4p+1)\times (4p+1)}.
   \end{equation}

At each iteration new columns are added to the matrices $R_L$ and $R_R$. To further
reduce the computational cost, we update the QR decomposition so as to
avoid a new factorization from scratch at each iteration.
   Assume that we have a current subspace of size $p$, so that
 $V_p = [v_1,\hdots, v_p]$ and we add another vector $v_{p+1}$ to the space. 
Thus, four new columns are added to $R_L$, which we add to the end of $R_L$, 
while keeping the same reordering for $R_R$. Adding a new column $v_{p+1}$ gives
   \begin{equation}
   \begin{bmatrix} R_L & v_{p+1} \end{bmatrix} = \begin{bmatrix} Q_1 R_1 & v_{p+1}\end{bmatrix} = \begin{bmatrix} Q_1 & q_2 \end{bmatrix} \begin{bmatrix} R_1 & R_{1,2} \\ 0 & R_{2,2} \end{bmatrix},
   \end{equation}
   where only two column vectors $q_2$ and $R_{1,2}$ and a scalar value $R_{2,2}$ have to be computed. We have
   \begin{equation}
   Q_1 R_{1,2} + q_2 R_{2,2} = v_{p+1}. \label{eq:newQR}
   \end{equation}
   Setting $R_{1,2} = Q_1^T v_{p+1}$ and constructing the vector
   \begin{equation}
   \hat{q}_2 = v - Q_1(Q_1^T v_{p+1}) = v_{p+1} - Q_1 R_{1,2},
   \end{equation}
   we can set $q_2 = \frac{\hat{q}_2}{\|\hat{q}_2\|}$ and thus $R_{2,2} = \|\hat{q}_2\|$. With this, Equation \eqref{eq:newQR} holds and the new column $q_2$ is orthogonal to $Q_1$,
   \begin{equation}
   Q_1^T q_2 = Q_1^T \frac{\hat{q}_2}{R_{2,2}} = \frac{Q_1^T v_{p+1} - Q_1^T Q_1 Q_1^T v_{p+1}}{R_{2,2}} = 0.
   \end{equation}
   
   Now we can completely avoid forming the full residuals, as with $R = R_L R_R^T = Q_1 R_1 R_R^T$ the desired Frobenius norm of the residual becomes
   \begin{equation}
   \|R\|_F = \sqrt{\mbox{trace}(R^T R)} = \sqrt{\mbox{trace}(R_R R_1^T Q_1^T  Q_1  R_1  R_R^T)} = \sqrt{\mbox{trace}(R_R R_1^T R_1 R_R^T)}.  \label{eq:QR_final}
   \end{equation}
   To compute the residual norms in this form only a very small matrix of size $n_T \times n_T$ has to be computed. Both of these residuals can be computed without forming the approximations $Y$ and $\Lambda$.
The use of the trace square root may lead to inaccuracy at the level of the machine epsilon square root. 
However, our convergence tolerance is in general significantly larger, say $10^{-4}$, hence these problems were
not encountered.
   
This computational procedure allows us to cheaply compute both $\|R_1\|_F$ and $\|R_2\|_F$, that is the absolute value of the residual norms.
Following the guidelines in \cite{Higham2002} for linear matrix equations, we also
monitor a scaled backward error norm of the two 
equations, that is 
   \begin{equation}
   \begin{aligned}
    \rho_3 &= \frac{\|R_1\|_F}{\tau(\|M_1\|Z_Y\|_F + \|K\|_F \|Z_{\Lambda}\|_F + \|M\|_F \|\hat{Y}\|_F) + \|M\|_F\|Z_Y\|_F\|C\|_F}, \\
   &+ \frac{\|R_2\|_F}{\tau(\|K\|_F \|Z_Y\|_F + \frac{1}{\beta}\|M\|_F\|Z_{\Lambda}\|_F + \|M\|_F\|Z_Y\|_F\|C\|_F}. \label{eq:resFull}
   \end{aligned}
   \end{equation}
which takes into account the data order of magnitude. Summarizing, we stop our iteration whenever
$$
\max \{ \|R_1\|_F, \|R_2\|_F, \rho_3 \} \le {\tt tol}
$$
where {\tt tol} is a prescribed tolerance.

\section{Numerical Results}\label{sec:num}
We now present the performance and flexibility of our method on multiple examples for different PDE-constrained optimization problems. The spatial FE-discretizations of the PDE operator were conducted using the deal.II framework \cite{deal2007} with Q1 finite elements and the method described in the previous section was implemented in \matlab \, R2018b. All experiments were run on a desktop computer with an \intel \, Core i7-4770 Quad-core processor running at 4 $\times$ 3400 MHz with 32 GB of RAM. 

We will show robustness with respect to the discretization sizes $n$ and $n_T$ as well as the control parameter $\beta$. Furthermore, we demonstrate the applicability of our method to multiple different setups such as partial observation of the desired state, boundary control and a different non-symmetric PDE-operator.

{Other low-rank solvers.} To emphasize the competitiveness of our new approach we compare it with another low-rank method aimed at the same problem class. The idea of exploiting the system structure to avoid forming a large system of equations and finding a low-rank way to compactly represent the solution is not an entirely new approach. Earlier, in \cite{stoll2014} the authors developed a technique rewriting the problem in a low-rank context and solving it with a preconditioned Krylov subspace solver, namely \minres\ introduced in \cite{paige1975}. We denote this solver as \lrminres\ (low-rank \minres) \footnote[1]{The \lrminres\ MATLAB code used for our comparison is available under \hyperlink{https://www.tu-chemnitz.de/mathematik/wire/pubs/PubCodes.zip}{https://www.tu-chemnitz.de/mathematik/wire/pubs/PubCodes.zip}.}.
	
	Here, the variables are represented in a low-rank matrix format rather than vectors, as $Y = W_Y V_Y^T$, $U = W_U V_U^T$ and $\Lambda = W_{\Lambda} V_{\Lambda}^T$ and the system of equations \eqref{eq:system1} - \eqref{eq:system3} is split into a matrix product as
		\begin{equation}
		\begin{bmatrix} \tau M_1 W_Y & \tau K^T W_{\Lambda} & M \Lambda \end{bmatrix} \begin{bmatrix} V_Y^T \\ V_{\Lambda}^T \\ V_{\Lambda}^T C \end{bmatrix} = \tau M_1 \hat{Y}
		\end{equation}
		for Equation \eqref{eq:system1} and respectively for Equations \eqref{eq:betau} and \eqref{eq:system3}. From here, the system of equations is solved with a tailored version of \minres. Opposed to the original version of \minres, the residual and all other upcoming variables are formulated as low-rank matrices rather than vectors and subsequently truncated to keep memory requirements low. Combined with exploiting the structure of the low-rank equation system in the upcoming matrix vector products and preconditioning, this method provides a memory efficient and fast alternative to established approaches for solving PDE constrained optimization problems. 
		
A bottleneck in this approach, however, is the setup of a preconditioner. To construct the necessary low-rank preconditioner, a matrix equation needs to be solved which significantly affects the overall performance. Therefore, our approach of directly starting from a matrix equation formulation proofs to be superior in many setups.
	
\subsection{Fully observed heat equation}
As a first example PDE we use the heat equation with a full observation, thus $\Omega_1 = \Omega$ and $M_1 = M$. This leads to a simplified version of Equation \eqref{eq:matrix_equation} as
\begin{equation}
A_1X + X C_1 + X I_0 + A_3 X D - Y_1 F_2^T = 0.
\end{equation}
First, we will use this simple setup to investigate the convergence behavior and our choice of stopping criteria by comparing the results to a solution obtained by a direct solver accurate up to machine precision. The constant-in-time desired state $\hat{Y}$ is displayed in Figure \ref{figure:desired_state}. This constant desired state leads to a low-rank representation of the right hand side of rank 1. The reference solution was obtained by solving the linear equation system of a small setup with MATLABs \textit{backslash} operator. We used a discretization of $n = 1089$ and $n_T = 100$.

\begin{figure}[htb]
	\centering
	\includegraphics[width=0.9\textwidth]{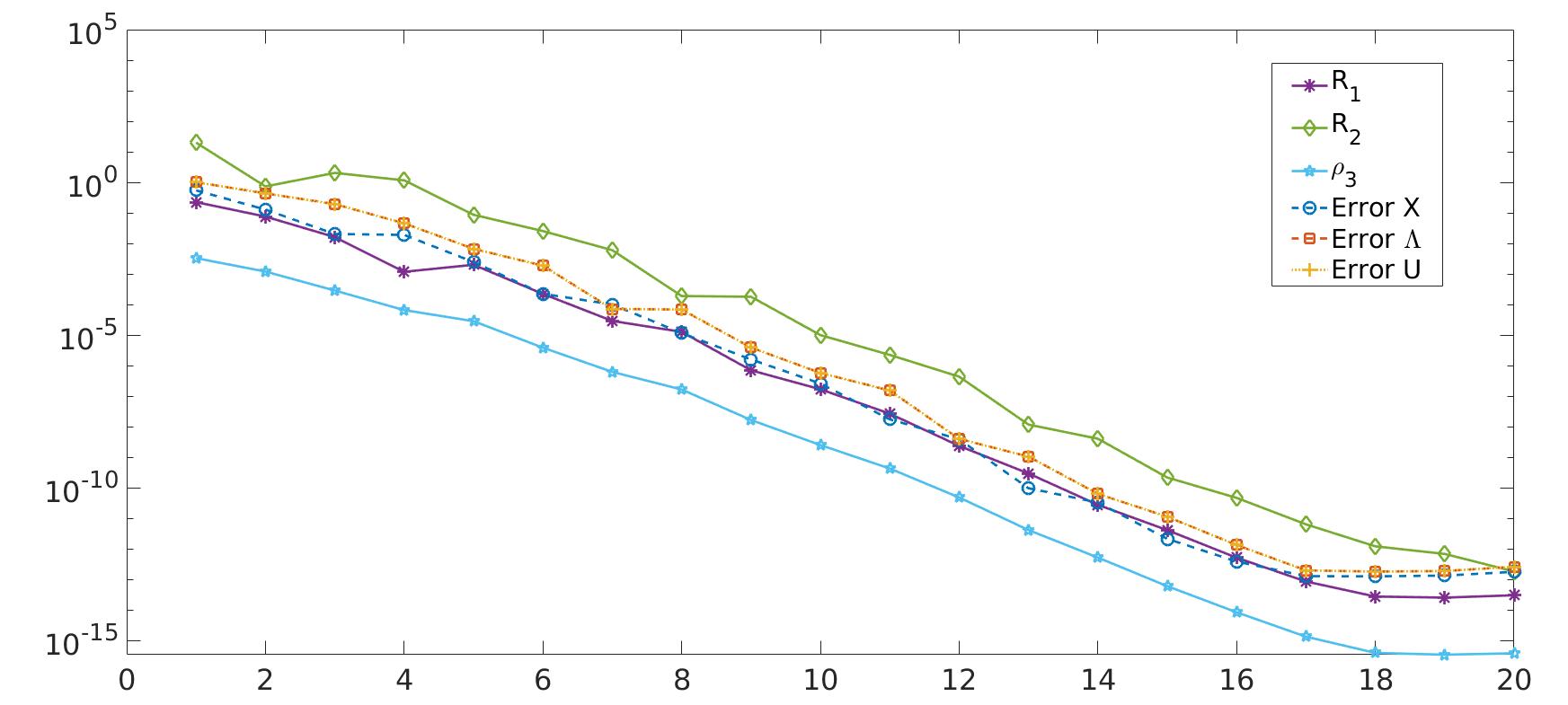}
	\caption{Observed stopping criteria and actual error progression.} \label{figure:convergence}
\end{figure}

The monitored quantities \eqref{eq:res1}, \eqref{eq:res2}, \eqref{eq:resFull} and the actual relative errors to the reference solution are displayed in Figure \ref{figure:convergence}. We see that the monitored quantities are a good estimation for the magnitude of the errors as $R_1$ stays close above  the actual errors. 
Inserting the direct solution into Equation \eqref{eq:matrix_equation_final} gives a residual of $5.11 \cdot 10^{-9}$. Thus, the different scaling of the matrix equation prohibits further accuracy gains once a high accuracy is reached. This effect is reflected in Figure \ref{figure:convergence} as the last iterations do not provide further accuracy gains.

\begin{example}\label{ex:constant-in-time}
{	\rm
We continue with the above simple example, i. e. the heat equation with full observation, $\Omega_1 = \Omega$ and $M_1 = M$ and a constant-in-time desired state.
We now investigate the performance of our method with respect to time and space discretization. We vary the number of time steps from 20 to 2500 and the number of spatial discretization nodes from $n=1024$ to $n=263169$ which roughly resembles up to a total of 78 million degrees of freedom. Here, again we fix the control parameter to $\beta = 10^{-4}$. As seen in Table~\ref{table:time_varying}, increasing the discretization size barely impacts the resulting subspace sizes. Additionally, the time needed to solve the optimization problems increases considerably slowly.
}
\end{example}

\begin{table}
	\centering
	\begin{tabular}{l l l l l l l l l l  r}
		\toprule
		$n$ & \multicolumn{2}{c}{1024} &  \multicolumn{2}{c}{4225} &  \multicolumn{2}{c}{16641} & \multicolumn{2}{c}{66049} &\multicolumn{2}{c}{263169} \\
		$n_T$ & $p$ & time(s) & $p$ & time(s) & $p$ & time(s) & $p$ & time(s) & $p$ & time(s)\\
		\midrule
20 & 8 & 0.09 & 11 & 0.10 & 11 & 0.26 & 11 & 1.02 & 12 & 5.82 \\ 
100 & 7 & 0.05 & 11 & 0.24 & 11 & 0.29 & 11 & 3.45 & 12 & 5.90 \\ 
500 & 6 & 0.05 & 11 & 0.23 & 11 & 0.42 & 11 & 1.16 & 13 & 6.61 \\ 
2500 & 6 & 0.19 & 11 & 0.91 & 11 & 1.09 & 10 & 1.61 & 15 & 9.52 \\
		\bottomrule
	\end{tabular}
\caption{Example \ref{ex:constant-in-time}
Subspace size $p$ and required time for different discretizations.} \label{table:time_varying}
\end{table}

\begin{table}[H]
	\centering
	\begin{tabular}{l l | r r r |r r r |r r r}
		\toprule
		$n$ &  & \multicolumn{3}{c}{1089} &  \multicolumn{3}{c}{4225} &  \multicolumn{3}{c}{16641} \\ 
		$\beta$ & method &  p & r & time(s) & p & r & time(s) &  p & r & time(s) \\ 
		\midrule
		$10^{-1}$ & \newmethod\  & 6 & 6& 0.09 & 5 & 5& 0.05 & 5 & 5& 0.14 \\ 
		& \lrminres\ & 20 & 17 & 4.00 & 22 & 21 & 7.49 & 22 & 21 & 31.59 \\ 
		$10^{-3}$ & \newmethod\ & 7 & 7& 0.07 &  7 & 7& 0.09 &  7 & 7& 0.22  \\ 
		& \lrminres\ &  11 & 7 & 1.80 & 14 & 11 & 5.46 &  11 & 10 & 24.75  \\ 
		$10^{-5}$ & \newmethod\ & 14 & 13& 0.19 &14 & 13& 0.23 &  14 & 13& 0.53 \\ 
		& \lrminres\ & 7 & 6 & 1.08 &  7 & 6 & 4.75 & - & - & - \\ 
		\hline 
	\end{tabular}
	\caption{Example~\ref{ex:comparisons}.
Comparison between the new \newmethod and another low-rank scheme \lrminres regarding 
required memory and time for $\hat{Y}$ with rank 1. \label{table:minres}}
\end{table}

\begin{example}\label{ex:comparisons}
{\rm
With the same data as in Example~\ref{ex:constant-in-time} we report on 
performance comparisons with the low-rank preconditioned \minres\ (\lrminres) proposed in \cite{stoll2014} and introduced above.
We fixed the number of time steps to $n_T = 100$ and computed solutions for different discretization sizes and control parameters $\beta$. The maximum discretization size here is $n=16641$.  
The results in Table~\ref{table:minres} reveal that our method's performance is superior regarding both required memory and CPU time. 
Here, our method is labeled as \newmethod. {The column denoted by $p$ states the 
subspace size for \newmethod, while for \lrminres\ it denotes the maximum rank per vector needed during the iterations of which up to 15 are required. The column $r$ states the rank of the final solution in both methods. Note that even though sometimes the ranks achieved by \lrminres\ are smaller the required memory to store the solution is still greater. This is because the \lrminres\ scheme needs to store a low-rank representation $U_1 U_2^T$ for each of the three variables as opposed to \newmethod\ where the same subspace $V_p$ is used for all variables. 
Additionally, during the iterations our method needs to store only one subspace of size $p$,
whereas \lrminres\ requires to store multiple vectors of size $p$. Also note, that our scheme} does not rely on the time consuming computation of a preconditioner as in \lrminres\ which requires the solution of a matrix equation that gets increasingly difficult for large spatial discretizations combined with small $\beta$. For the last entry of Table \ref{table:minres} \lrminres\ did not reach convergence.

 With the same settings as before, we now raise the rank of the desired state $\hat{Y}$ to 6 -- leading to a larger rank matrix $F_1 F_2^T$ -- and compare both methods once again. Table \ref{table:minres2} displays the results as before. Again, our method successfully converges with a small subspace sizes 
within a very short amount of time, whereas \lrminres\ 
was considerably slower especially for larger problem sizes and did not converge for the largest discretization with $\beta$ being very small.

\begin{table}[htb]
	\centering
	\begin{tabular}{l l |r r r  |r r r  |r r r }
		\toprule
		$n$ &  & \multicolumn{3}{c}{1089} &  \multicolumn{3}{c}{4225} &  \multicolumn{3}{c}{16641} \\ 
		$\beta$ & method  & $p$ & r & time(s) & $p$ & r & time(s) & $p$ & r & time(s) \\ 
		\midrule
$10^{-1}$ & \newmethod\ &  21 & 21& 0.38 & 21 & 21& 0.47 & 19 & 19& 0.72  \\ 
& \lrminres\ &  29 & 22 & 20.77 &  27 & 25 & 13.38 & 29 & 29 & 23.63  \\ 
$10^{-3}$ & \newmethod\ &  31 & 31& 1.08 &  30 & 30& 1.09 & 30 & 29& 1.77  \\ 
& \lrminres\ &  15 & 11 & 2.22 &  17 & 11 & 11.19 & 19 & 15 & 31.05  \\ 
$10^{-5}$ & \newmethod\ &  49 & 43& 3.75 &  49 & 44& 3.97 & 51 & 41& 5.88  \\ 
& \lrminres\ & 11 & 8 & 1.49 &  11 & 9 & 6.14 & - & - & -  \\ 
	\end{tabular}
	\caption{Example~\ref{ex:comparisons}. Comparison of  \newmethod and \lrminres , regarding required memory and CPU time for rank($\hat{Y}$)=6. \label{table:minres2}}
\end{table}

We are also interested in the memory savings our method provides. 
Therefore, we monitor the memory consumption of \newmethod, \lrminres\ and a full rank \minres\ method for the 
same set-ups displayed in Table~\ref{table:minres2}. For \newmethod\ we monitor the memory needed to 
construct the subspace, the reduced solution and the setup of the reduced equation system \eqref{eq:lineq}. For \lrminres\ we monitor the memory for all vectors that are used. For comparison we additionally report the memory consumed by a standard \minres method only to store the required vectors, not taking into account the system matrix and preconditioner. 
The results are shown in Table~\ref{table:memory}. The quantities refer to the overall memory 
consumption  in megabytes during the process of solving the equation system with the respective method. We see that even with constructing the reduced equation system \ref{eq:lineq} the memory requirements stay well below those of a not low-rank approach.}
\end{example}

\begin{table}[htb]
	\centering
	\begin{tabular}{l l | r r r } 
		\toprule
		& &  \multicolumn{3}{c}{Memory (MB)} \\ 
		n & $\beta$ & \newmethod\ & \lrminres\ & {\sc minres} \\
		\midrule
		1089 & 1e-1  & 4.46 & 8.11 & 15.39\\ 
		& 1e-3  & 8.06 & 6.14 &  15.39\\ 
		& 1e-5  & 17.58 & 4.54 &  15.39\\ 
		4225 & 1e-1  & 8.58 & 19.21 & 60.46\\ 
		& 1e-3  & 12.57 & 19.35 & 60.46\\ 
		& 1e-5  & 24.36 & 16.15 & 60.46\\  
		16641 & 1e-1  & 23.50 & 71.74 & 241.84\\ 
		& 1e-3  & 31.86 & 80.56 & 241.84 \\ 
		& 1e-5  & 51.91 & 205.85 & 241.84\\ 
	\end{tabular}
	\caption{Example~\ref{ex:comparisons}.
Memory requirements of all considered methods. \label{table:memory}}
\end{table}

\subsection{Partial observation on heat equation}

Until now we only investigated the case where $M = M_1$ and therefore the matrix $A_2$ being the identity. One highly interesting and challenging problem is the optimization under partial observation, meaning the desired state is only of interest on a partial domain $\Omega_1 \subseteq \Omega$. This leads to a singular matrix $M_1$, which further increases the difficulty of the PDE-constrained optimization problem.When we set up the rational Krylov subspace with only $A_1$ for this setting, we do not reach convergence. Therefore, we add up to two new vectors in each iteration as outlined in Section \ref{section:comp_specifications1}. 

\begin{figure}[htb]
	\centering
	\begin{subfigure}{0.32\textwidth}
		\includegraphics[width=0.9\textwidth]{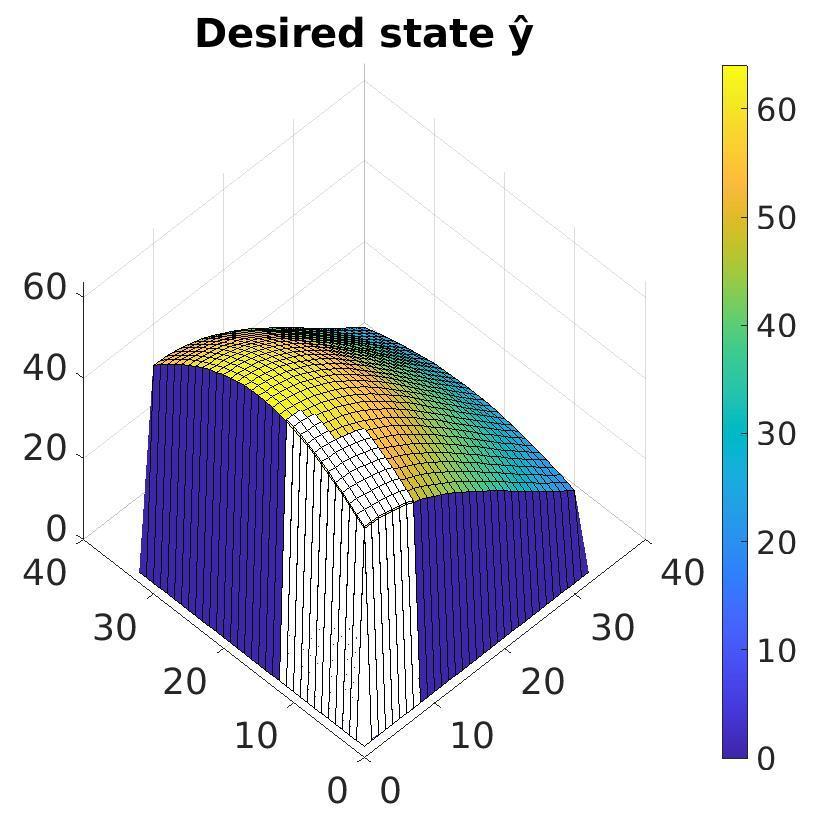}
		\caption{Partially observed desired state} \label{figure:desired_state}
	\end{subfigure}
	\begin{subfigure}{0.32\textwidth}
		\includegraphics[width=0.9\textwidth]{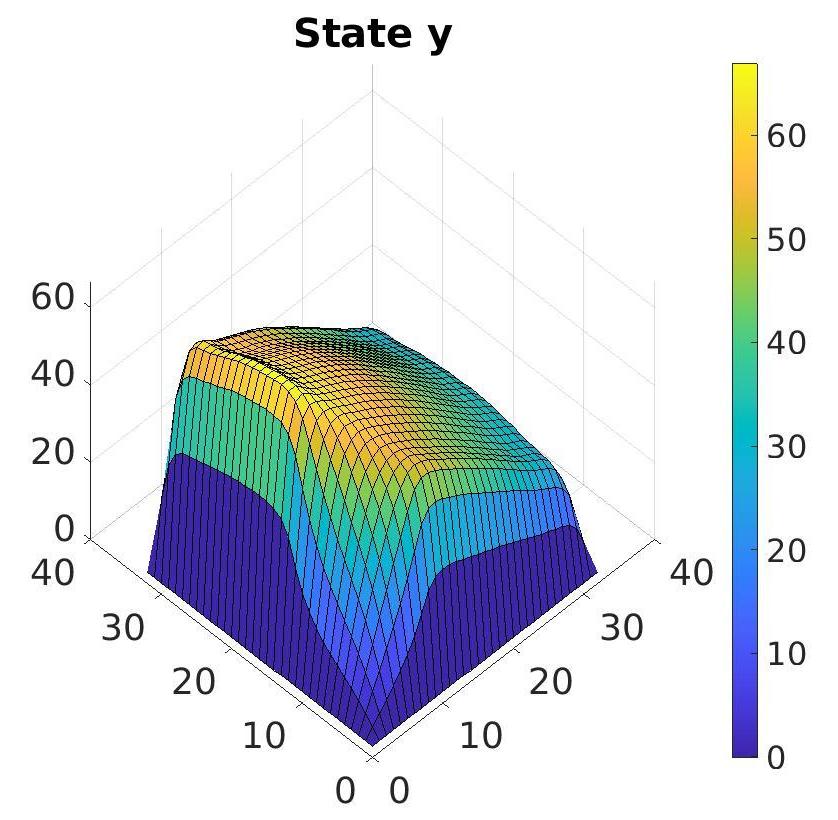}
		\caption{State solution for partial observation} \label{figure:partial_state}
	\end{subfigure}
	\begin{subfigure}{0.32\textwidth}
		\includegraphics[width=0.9\textwidth]{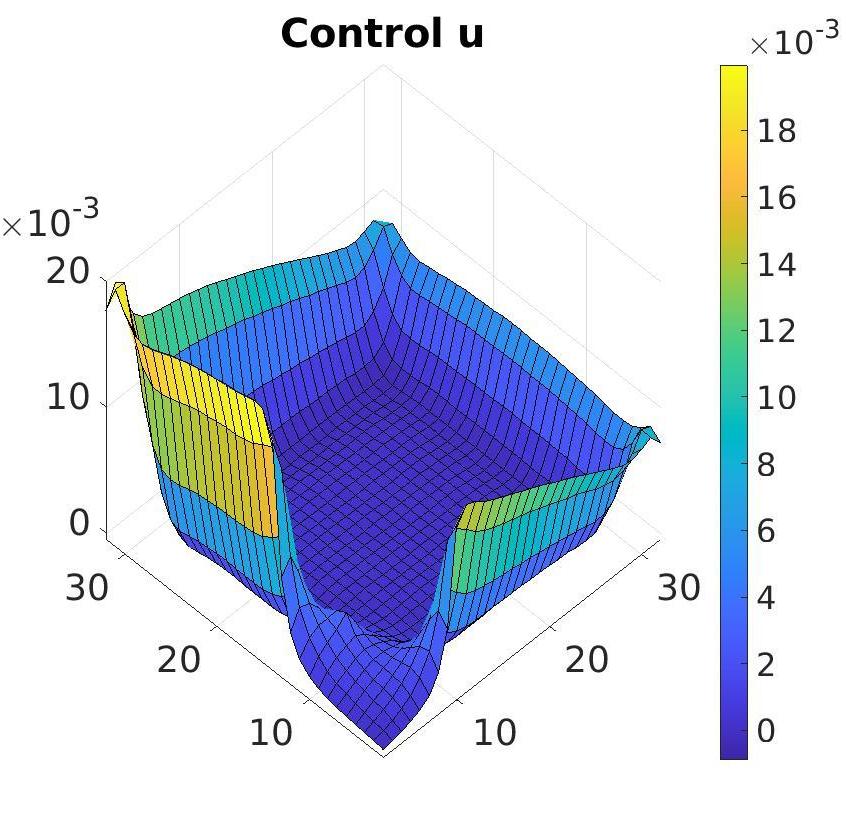}
		\caption{Control solution for partial observation} \label{figure:partial_control}
	\end{subfigure}
\caption{Example \ref{ex:partial_obs}. Solution at a fixed time instant for a partial observation. \label{figure:partial}}
\end{figure}

\begin{example}\label{ex:partial_obs}
{\rm
We solved the problem with control parameter $\beta = 10^{-4}$ and 100 time steps on a grid of {1089} spatial DoFs. 
In Figure~\ref{figure:desired_state} we see the desired state for this problem and the white 
area is an example of $100$ non-observed nodes, {roughly 10 \%}. Hence, 
the corresponding entries in $M_1$ are set to 0. Figure~\ref{figure:partial_state} and \ref{figure:partial_control} show the resulting state and control respectively for a fixed time step $t=25$.
Table~\ref{table:partial} shows the results when
increasing the number of unobserved nodes from $n_0 = 0$ to $n_0 = 900$ which is about 90\% of the nodes. 
Here, the constructed subspaces are not optimal. 
Thus, we make use of the truncation modification outlined in Chapter \ref{section:comp_specifications1}.  We see that the time and iterations needed to solve this more challenging problem are higher than for the fully observed case in the table's first row. 
For the non-truncated case with full $V_p$, the columns denoted by $p$ for the subspace size and $r$ for the solution rank show a discrepancy which disappears with stricter truncation. 
 When truncating with a tolerance of $10^{-10}$ both values coincide for all cases. 
In some cases applying this truncation greatly increases the number of iterations needed to find a solution but the end result is a higher memory reduction. Thus, this option can be toggled to better reflect the desired behavior.

 Unfortunately, the \lrminres\ scheme from \cite{stoll2014} is not adapted for the case of $M \neq M_1$. Hence, we have no other low-rank method to compare the performance of our method  with.
}
\end{example}

\begin{table}[htb]
	\begin{tabular}{ l r r r r r r r r r r r r }
		\toprule
		& \multicolumn{4}{c}{full $V_p$} & \multicolumn{4}{c}{truncated $V_p$, $10^{-12}$} & \multicolumn{4}{c}{truncated $V_p$, $10^{-10}$}  \\
		$n_0$ & time(s) & iters & $p$ & $r$ & time(s) & iters & $p$ & $r$ &  time(s) & iters & $p$ & $r$ \\
		\midrule
0 & 0.06 & 7 & 7 & 7& 0.05 & 7 & 7 & 7& 0.06 & 7 & 7 & 7\\ 
100 & 0.67 & 19 & 31 & 29& 0.67 & 19 & 30 & 29& 0.69 & 19 & 23 & 23\\ 
300 & 1.46 & 26 & 42 & 34& 1.41 & 26 & 35 & 34& 1.20 & 28 & 26 & 26\\ 
500 & 2.55 & 32 & 51 & 37& 2.52 & 34 & 38 & 36& 24.41 & 142 & 29 & 29\\ 
700 & 1.50 & 26 & 42 & 37& 1.59 & 27 & 38 & 36& 2.46 & 42 & 29 & 29\\ 
900 & 0.82 & 21 & 34 & 34& 0.84 & 21 & 33 & 33& 0.75 & 21 & 25 & 25\\ 
		\bottomrule
	\end{tabular}
	\caption{Example \ref{ex:partial_obs}. Results for different levels of partial observation. \label{table:partial}}
\end{table}

\subsection{Boundary Control}
Another setup of high interest is having a non-distributed control. Here, we take a boundary control problem as an example. This leads to the control $u$ not being distributed on $\Omega$ but only on a subdomain $\Omega_b \subset \Omega$. Therefore, we get a 
smaller mass matrix $M_b \in \R^{n_b \times n_b}$ associated with $U$ present on $n_b < n$ spatial nodes,
 and $N \in \R^{n \times n_b}$ is now rectangular (tall). 
Therefore, we have a modification of Equation \eqref{eq:matrix_equation_final} as
\begin{equation}
A_1 X +  X C_1 + X I_0 + A_3 X D - F_1 F_2^T = 0
\end{equation}
with $A_1 = M^{-1}K$ and $A_3 = M^{-1}N M_b^{-1} N^T$. Here, the subspace we compute is derived from $A_1$ and $A_3$ as in \eqref{eq:subspace}.

\begin{example}\label{ex:bc}
{\rm
 As before, we use $n_T = 100$ time steps and a convergence threshold of $10^{-4}$ for this setup. Results with a constant in time desired state for different levels of spatial discretization and 
different values for $\beta$ are displayed in Table \ref{table:minresbc}. \newmethod\ produces robust results with respect to the discretization size as we can only see a small increase in the ranks throughout discretization. 
Additionally, the results are acquired within a short amount of time even for small $\beta$. For 
larger discretization sizes the required time grows only slightly 
opposed to \lrminres\ where time requirements increase substantially for larger discretization sizes.
}
\end{example}

\begin{table}[H]
	\centering
	\begin{tabular}{l l | r r r | r r r | r r r} 
		\toprule
		$\quad n$ &  & \multicolumn{3}{c}{1089} &  \multicolumn{3}{c}{4225} &  \multicolumn{3}{c}{16641} \\ 
		$\beta$ & method &  p & rank & time(s) & p & rank & time(s) & p & rank & time(s) \\ 
		\midrule
$10^{-1}$ & \newmethod\ & 72 & 47& 3.92 &  80 & 47& 7.33 &  80 & 47& 16.35  \\ 
& \lrminres\ & 29 & 15 & 29.47 &  28 & 16 & 113.99 &  28 & 16 & 496.70  \\ 
$10^{-3}$ & \newmethod\ &  64 & 46& 3.68  & 68 & 46& 5.08 &  74 & 46& 14.35  \\ 
& \lrminres\ &  25 & 12 & 12.70  & 24 & 14 & 35.53 & 24 & 14 & 158.82  \\ 
$10^{-5}$ & \newmethod\ &  104 & 45& 12.36  & 116 & 44& 18.95 &  107 & 44& 28.48  \\ 
& \lrminres\ &  21 & 16 & 34.20 &  21 & 16 & 142.83  & 20 & 16 & 694.93  \\ 
			\bottomrule
	\end{tabular}
	\caption{Example \ref{ex:bc}. Comparisons between the new \newmethod\  and \lrminres\ 
regarding required memory and time for a boundary control problem.} \label{table:minresbc}
\end{table}

\subsection{Non-symmetric convection-diffusion PDE}

Until now we only investigated the performance of our method under the constraint of a heat equation. But for more challenging PDEs our method works just as well, e.g. the convection-diffusion equation,
\begin{equation}
 \dot{y} = \eps \nabla^2 y - v \cdot \nabla y + u,
\end{equation}
with diffusion coefficient $\eps$ and velocity field $v$. 
When $\eps \ll 1$, the equation is convection-dominated and solving it is a challenging task. 
 As the stiffness matrix $K$ of the resulting equation system is non-symmetric, we have to adjust the matrix equation accordingly. 
Thus, we get another modification of Equation \eqref{eq:matrix_equation_final} for this setup as
\begin{align}
A_1 X \begin{bmatrix} \mathcal{I} & 0  \\ 0 & 0 \end{bmatrix} + M^{-1} K^T X \begin{bmatrix} 0 & 0  \\ 0 & \mathcal{I} \end{bmatrix} + X C_1 + A_2 X I_0 + A_3 X D - F_1 F_2^T = 0.
\end{align}

\begin{figure}[htb]
	\centering 
	\begin{subfigure}[b]{0.32\textwidth}
		\includegraphics[width=1\textwidth]{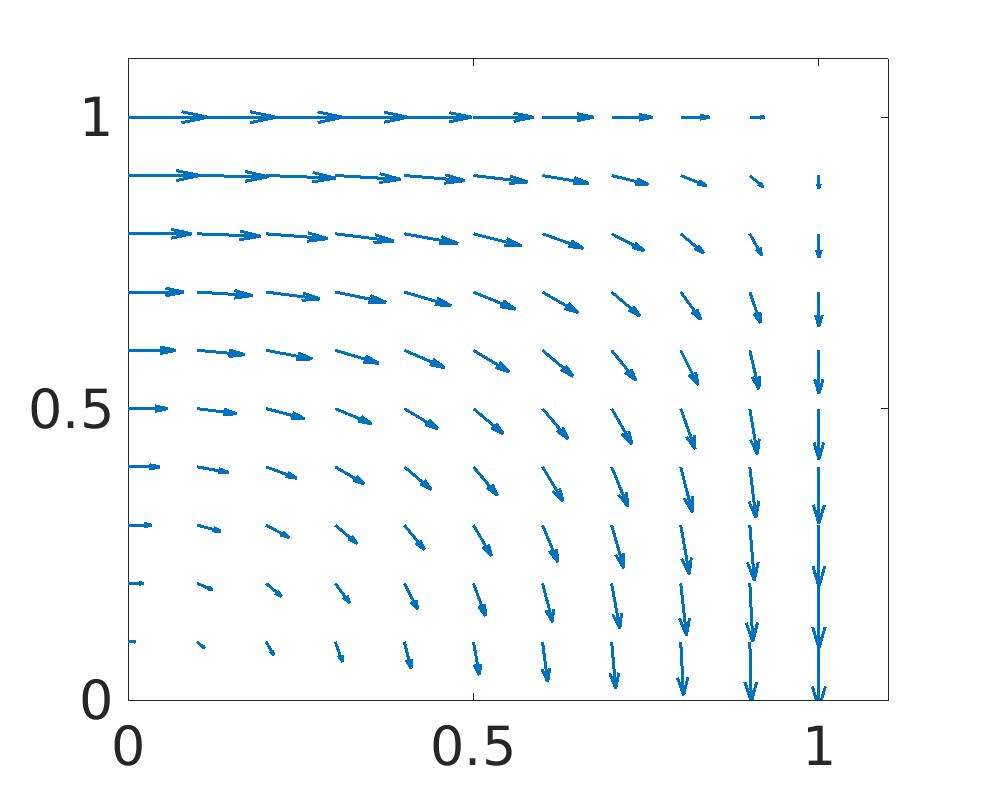}
		\caption{Velocity field $v$} \label{figure:wind}
	\end{subfigure}
	\begin{subfigure}[b]{0.32\textwidth}
		\includegraphics[width=0.9\textwidth]{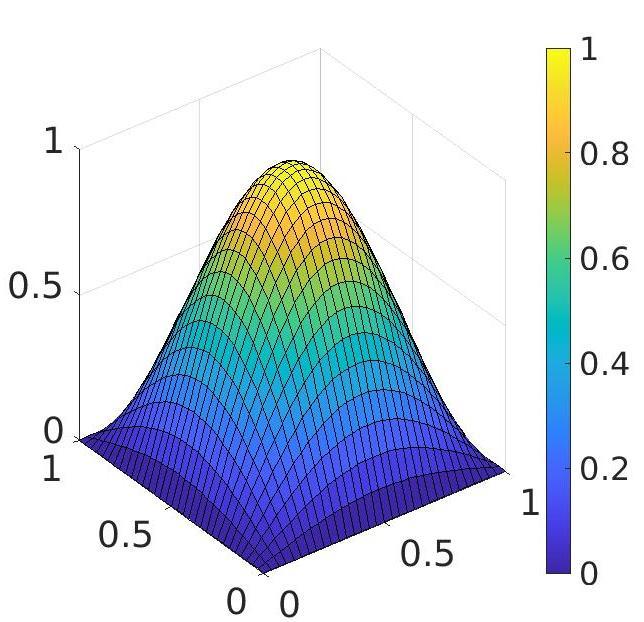}
		\caption{Desired state} \label{figure:CD_state}
	\end{subfigure}
	\begin{subfigure}[b]{0.32\textwidth}
		\includegraphics[width=0.9\textwidth]{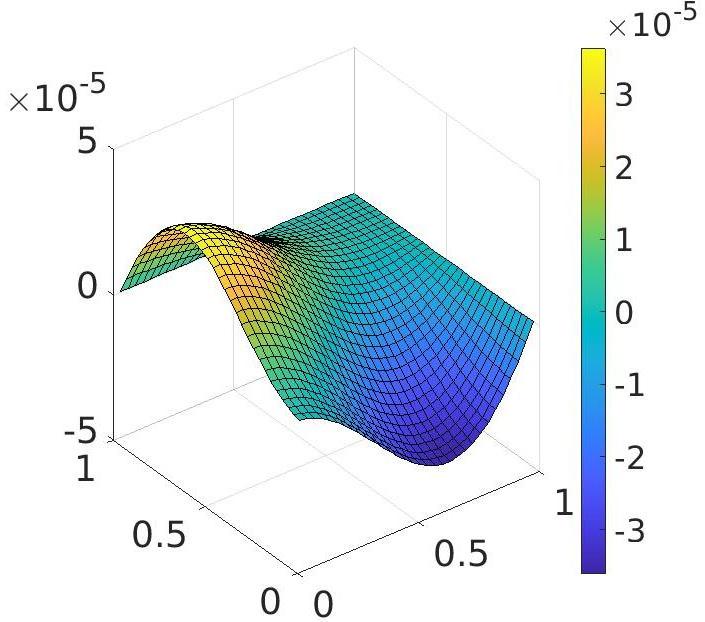}
		\caption{Control solution} \label{figure:CD_control}
	\end{subfigure}
	\caption{Example~\ref{ex:cd}. Solution at a fixed time instant for convection-diffusion problem. \label{figure:CD}}
\end{figure}

\begin{example}\label{ex:cd}
{\rm
We consider example 3.1.4 from \cite{elman2005} with a recirculating 
wind $v(x,y) = (2y(1-x^2) , -2x(1-y^2))$ as the underlying model. The velocity field $v$ is displayed in Figure~\ref{figure:wind}.
 In Figures~\ref{figure:CD_state}-\ref{figure:CD_control} the desired state and a snapshot of a control solution for $\beta = 10^{-5}$ are displayed for $n=4225$ and $n_T = 100$. The \newmethod\ method produces reliable results throughout a range of different values for $\eps$ and $\beta$ as reported in Table~\ref{table:CD}. 
Even very small values in $\eps$ do not pose a problem to the solver. 
}
\end{example}

\begin{table}[htb]
	\centering
	\begin{tabular}{l l |l r r |l l r|l r r}
		\toprule
		$n$ &  & \multicolumn{3}{c}{4225} &  \multicolumn{3}{c}{16641} &  \multicolumn{3}{c}{66049} \\
		$\beta$ & $\eps$ & mem & $p$ & time & mem &$p$  & time & mem &$p$ & time \\ 
\midrule
 $10^{-1}$  & $ 1$ & 0.03 & 5 & 0.25& 0.01 & 1 & 0.43& 0.01 & 1 & 4.61 \\ 
 & $ 10^{-1}$ & 0.05 & 9 & 0.27& 0.07 & 13 & 1.51& 0.05 & 9 & 7.87 \\ 
 & $ 10^{-2}$ & 0.15 & 29 & 1.42& 0.10 & 19 & 2.56& 0.10 & 19 & 15.40 \\ 
 & $ 10^{-3}$ & 0.23 & 43 & 3.04& 0.09 & 17 & 2.96& 0.09 & 17 & 17.16 \\ 
 \midrule
 $10^{-3}$  & $ 1$ & 0.06 & 11 & 0.33& 0.05 & 9 & 1.10& 0.05 & 9 & 9.07 \\ 
 & $ 10^{-1}$ & 0.09 & 17 & 0.51& 0.09 & 17 & 2.17& 0.05 & 9 & 7.94 \\ 
 & $ 10^{-2}$ & 0.08 & 15 & 0.55& 0.10 & 19 & 2.71& 0.07 & 13 & 10.60 \\ 
 & $ 10^{-3}$ & 0.34 & 65 & 8.29& 0.07 & 13 & 2.37& 0.09 & 17 & 16.97 \\ 
 \midrule
 $10^{-5}$  & $ 1$ & 0.10 & 19 & 0.70& 0.08 & 15 & 2.08& 0.08 & 15 & 13.53 \\ 
 & $ 10^{-1}$ & 0.02 & 3 & 0.08& 0.02 & 3 & 0.38& 0.16 & 31 & 26.91 \\ 
 & $ 10^{-2}$ & 0.04 & 7 & 0.18& 0.02 & 3 & 0.48& 0.02 & 3 & 3.01 \\ 
 & $ 10^{-3}$ & 0.08 & 15 & 0.56& 0.02 & 3  & 1.06& 0.02 & 3 & 6.10 \\ 
		 \bottomrule 
		\end{tabular}
	\caption{Example~\ref{ex:cd}. Results for different diffusion 
coefficients $\eps$ for the convection-diffusion problem with different 
parameters $\beta$ and different spatial discretizations. \label{table:CD}}
\end{table}

\section{Conclusion and Outlook}\label{sec:conc}
We proposed a new scheme to solve  a class of
large-scale PDE-constrained optimization problems in a low-rank format. This method relied on the reformulation of the KKT saddle point system into a matrix equation, which was subsequently projected onto a low dimensional space 
 generated with rational Krylov type iterations. 
We showed that the method's convergence is robust with respect to different discretizations and parameters. Furthermore we demonstrated higher memory and time savings compared to an established low-rank scheme. Additionally, the \newmethod\ is very flexible with respect to different constraints such as non-symmetric PDE operators or partial state observation.

In the future, we plan on further investigating the subspace choice and the performance of our 
scheme for other challenging setups. Further improvements 
are expected from realizing a truncation or restarting mechanism in cases where the subspaces become unexpectedly large.

\section*{Acknowledgments}
The work of the first and third authors was supported by the German ScienceFoundation (DFG) through grant 1742243256 - TRR 9.
The second author is a member of Indam-GNCS. Its support is gratefully acknowledged.
This work was also supported by 
the DAAD-MIUR Mobility Program 2018 ``Optimization and low-rank solvers for isogeometric analysis (34876)'' between
TU-Chemnitz and the Universit\`a di Bologna.

\bibliographystyle{siamplain}
\bibliography{references}

\end{document}